\documentclass{article}
\usepackage{amsthm}
\usepackage[utf8]{inputenc}
\usepackage[T1]{fontenc}
\usepackage[centertags]{amsmath}
\usepackage{amsfonts}
\usepackage{graphicx}
\usepackage{amsmath}
\usepackage{amssymb}
\usepackage{latexsym}
\usepackage[mathscr]{euscript}
\usepackage{tikz,color,makeidx} 
\usetikzlibrary{calc,matrix,arrows,chains,positioning,scopes,decorations.pathmorphing,backgrounds,fit}
\numberwithin{equation}{section}
\newtheorem{teo}{Theorem}[section]
\newtheorem{cor}{Corollary}[section]
\newtheorem{pro}{Proposition}[section]
\newtheorem{lem}{Lemma}[section]
\newtheorem{es}{\textbf{Example}}[section]
\newtheorem{defi}{\mbox{\textbf{Definition}}}[section]
\newtheorem{rem}{\mbox{\textbf{Remark}}}[section]
\newtheorem{notation}{\textit{Notation}}[section]

\newcommand{\bdfn}{\begin{defi} \begin{rm}}
\newcommand{\edfn}{\end{rm} \end{defi}}
\newcommand{\bthm}{\begin{teo}}
\newcommand{\ethm}{\end{teo}}
\newcommand{\bprop}{\begin{pro}}
\newcommand{\eprop}{\end{pro}}
\newcommand{\bcor}{\begin{cor}}
\newcommand{\ecor}{\end{cor}}
\newcommand{\blem}{\begin{lem}}
\newcommand{\elem}{\end{lem}}
\newcommand{\bfact}{\begin{rem} \begin{rm}}
\newcommand{\efact}{\end{rm} \end{rem}}
\newcommand{\bex}{\begin{es} \begin{rm}}
\newcommand{\eex}{ \end{rm} \end{es}}
\newcommand{\ten}{\otimes}
\newcommand{\quot}[2]{{\raisebox{.2em}{$#1$}\left/\raisebox{-.2em}{$#2$}\right.}}
\newcommand{\vsp}[1]{\vspace*{#1cm}}
\newcommand{\bnot}{\begin{notation} \begin{rm}}
\newcommand{\enot}{\end{rm} \end{notation}}

\tikzset{node distance=2cm, auto}

\title{\mbox{Stochastic independence for probability MV-algebras}}

\author{Serafina Lapenta\\ 
 {\small Department of Mathematics, Computer Science and Economics,}\\ 
{\small University of Basilicata,
 Viale dell'Ateneo Lucano, 10 Potenza, Italy}\\
 {\small serafina.lapenta@unibas.it}
\and 
Ioana Leu\c stean \\
{\small Department of Computer Science,} \\
{\small Faculty of Mathematics and Computer Science, University of Bucharest,}\\
{\small Academiei nr.14, sector 1, C.P. 010014,  Bucharest, Romania}\\   
{\small ioana@fmi.unibuc.ro}
}

\date{}

\begin{document}
\maketitle

\begin{abstract}
We prove that any   MV-algebra has a faithful state can be embedded in an \em{f}MV-algebra of integrable functions. As consequence,  we prove H\"{o}lder's inequality and Hausdorff moment problem for MV-algebras with product and  we  propose a solution for the stochastic independence of probability MV-algebras.
\end{abstract}

\section*{Introduction}

{\em MV-algebras} were defined by Chang \cite{Cha1} and they stand to \L ukasiewicz $\infty$-valued logic as boolean algebras stand to classical logic. The theory of MV-algebras was highlighted by
Mundici's categorical equivalence between MV-algebras and abelian
lattice-ordered groups with strong unit ($\ell u$-groups) \cite{Mun1}. The twofold nature of MV-algebras,  generalizations of boolean algebras and unit intervals of  $\ell u$-groups,  is also reflected by their probability theory: {\em the (finite-additive) states} are in one-to-one
correspondence with normalized states on $\ell u$-groups, while {\em probability MV-algebras} are the main
ingredient of the extension of  Carath\'{e}odory  boolean algebraic probability theory to many-valued events.

 A probability MV-algebra \cite{MunRie} is  a
pair $(A,s)$, where $A$ is a $\sigma$-complete MV-algebra  and $s$
is a  $\sigma$-continuous faithful state. Riec\v{a}n and Mundici  propose in \cite{MunRie}  a list of open problems, we recall the fifth one:

{\em"[...] Assuming  $M$ and $N$ to be probability MV-algebras, generalize the classical
theory of "stochastically independent" $\sigma$-subalgebras as defined
in Fremlin's treatise [Measure Theory, 325L]."}

\noindent In \cite{LeuStateC} the author investigates this problem for MV-algebras endowed with finite-additive states, but no solution is given for probability MV-algebras, as defined in \cite{MunRie}.  In the present paper, Theorem \ref{teo:unpropInd} presents a possible solution and we notice that an analogue result in Fremlin's treatise \cite{Fre} is [253F]. 

An important result in our approach is Theorem \ref{teo:embMeasure} which, combining the results from \cite{LeuStateC} and \cite{LeuSCRMV}, prove that any MV-algebra that has a faithful state can be embedded in an \textit{f}MV-algebra of integrable functions in which the state is represented by the integral. The representation for states is actually the Kroupa-Panti teorem \label{teo:KP},  but we make the context more precise. The representation of the algebraic structure is crucial for our development and  is based  on  Theorem \ref{teo:equivL1}, which is similar with Kakutani's representation for abstract $L$-spaces  \cite{Kak}. The
 {\em  {\textit f}MV-algebras} are defined in \cite{LL1} any they are MV-algebras endowed with both an internal product and a scalar product, with scalars from $[0,1]$. By an extension of Mundici's equivalence, they are categorically equivalent with unital \textit{f}-algebras \cite{BP}.  As direct consequences of  Theorem \ref{teo:embMeasure}, in Section \label{section:emb} we prove  {\em H\"{o}lder's inequality} and {\em Hausdorff's moment problem}  for {\em  PMV-algebras}, i.e. MV-algebras  endowed with an internal  product \cite{DiND},   and for {\textit f}MV-algebras.

Section \ref{section:ind} is focused on  the problem of {\em stochastic independence}. The main idea is the following: given two
probability MV-algebras  we  embed them in corresponding algebras of integrable functions, which allow us to apply the results 
from  \cite{Fre}.  Our final result can be stated as follows:

{\em Given   $(A, s_A)$ and $(B, s_B)$  two probability MV-algebras, there exists a probability MV-algebra  $(T,s_T)$  and a bilinear  function $\beta: A\times B\rightarrow T$ such that  $s_T(\lambda(a,b))=s_A(a)\cdot s_B(b)$, for any $a\in A$ and $b\in B$.}\\
By  Theorem \ref{teo:unpropInd}, the probability MV-algebra $(T,s_T)$ satisfy an universal property which, however, does not characterize it up to isomorphism.

\section{Preliminaries}

\subsection{Algebraic structures}
\noindent An {\em MV-algebra} is an algebraic structure $(A, \oplus, \ ^*, 0)$, where $(A, \oplus ,0)$ is a commutative monoid, $\ ^*$ is an involution and the relation $(a^*\oplus b)^*\oplus b=( b^*\oplus a)^*\oplus a$ is satisfied for any $a,b\in A$ \cite{Cha1, CDM, MunBook}. The variety of MV-algebras is generated by $([0,1],\oplus, ^*,0)$ where $a\oplus b=min(a+b,1)$ and $a^*=1-a$ for any $a,b\in [0,1]$. The category of MV-algebras is denoted $\mathbf{MV}$.

One  also defines  the constant $1=0^*$, the operation $a\odot b = (a^* \oplus b^*)^*$ and the distance function 
$d(a,b)=(a\odot b^*)\oplus (b\odot a^*)$ for any $a,b\in A$. 
Setting $a\le b$ if and only if $a^*\oplus b=0$, then $(A,\leq, 0,1)$ is a bounded distributive lattice such that 
$a\vee b=(a^*\oplus b)^*\oplus b$ and $a\wedge b=(a^*\vee b^*)^*$ for any $a,b\in A$.  An MV-algebra $A$  is  $\sigma$-complete (Dedekind-MacNeille complete) if its lattice reduct is a $\sigma$-complete (Dedekind-MacNeille complete) lattice.

If $A$ is an MV-algebra we define a partial operation $+$ as follows: for any $a,b\in A$, $a+b$ is defined if and only if $a\leq b^*$ and, in this case, $a+b=A\oplus b$. This operation is cancellative and  any MV-algebra $A$ satisfies the Riesz decomposition property \cite[Section 2.9]{LDHandbook}. Throughout the paper we use the following notation:
\vspace*{-0.2cm}
\begin{center}
$na=\underbrace{a+\cdots+ a}_n$ and $n_\oplus a=\underbrace{a\oplus\cdots\oplus a}_n$  
\end{center}
\vspace*{-0.2cm}
where $a\in A$ and $n\geq 1$ is a natural number.

If $A$ and $B$ are MV-algebras then a function $\omega:A\to B$ is  {\em linear}  if $f(a+b)=f(a)+f(b)$ whenever $a\leq b^*$. Bilinear functions are defined as usual. A bimorphism is a bilinear function that is $\vee$-preserving and $\wedge$-preserving in each component. We refer to \cite{FloLeu} for  basic results on linear functions.

An {\em ideal} in  $A$ is a lower subset $I$ that contains $0$ and it is closed to $\oplus$. A maximal ideal is an ideal that is maximal in the set of all ideals ordered by set-theoretic inclusion. A {\em semisimple MV-algebra} is an MV-algebra in which the intersection of all maximal ideals is $\{ 0 \}$. Equivalently, an MV-algebra $A$ is semisimple if and only if there exists a compact  Hausdorff space $X$ such that  $A$ can be embedded in the MV-algebra $C(X)=\{f:X\to [0,1]\mid f \mbox{ continuous}\}$ with pointwise operations \cite[Corollary 3.6.8]{CDM}.

A {\em $\ell$-group} is an abelian group that is also a lattice such that any group translation is isotone.   If $G$ is an $\ell$-group, an element $u \in G$  is a \emph{\textup{(}strong order\textup{)} unit} if $u\geq 0$, and for every $x\in G$ there is a natural number $n\geq 1$  such that $nu \geq |x|$. An $\ell$-group is {\em unital} if it is endowed with a distinguished unit.  The category of unital $\ell$-groups and their unital homomorphisms is denoted $\mathbf{auG}$. 

If $(G,u)$ is a  unital $\ell$-group we denote $[0,u]=\{x\in G| 0\leq x\leq u\}$ and  we define
\vsp{-0.2}
\begin{center}
$x\oplus y=(x+y)\wedge u$  and  $\neg x=u-x$ for any $x, y\in [0,u]$.
\end{center}
\vsp{-0.2}
Then $[0,u]_G=([0,u],\oplus,\neg,0)$ is an MV-algebra. For any MV-algebra $A$ there exists an $\ell u$-group $(G,u)$ such that $A\simeq [0,1]_G$. Moreover, the following property  holds: for any $x\geq 0$ in $G$ there exist a natural number $n\geq 1$ and $a_1,\ldots, a_n\in A$ such that $x=a_1+\cdots+a_n$.

 It is possible to define a functor $\Gamma: \mathbf{auG}\rightarrow \mathbf{MV}$ by 
\vsp{-0.2}
\begin{center}
$\Gamma (G,u)=[0,u]_G $ and $\Gamma(h)= h|_{[0,u]_G}.$ 
\end{center}
\vsp{-0.2}
where $(G,u)$ is an $\ell u$-group and $h$ is unital homomorphism.  In \cite{Mun1} it is proved that $\Gamma$ establishes a categorical equivalence between the categories $\mathbf{auG}$ and $\mathbf{MV}$. In addition, an MV-algebra  $A$ is semisimple if and only if the corresponding $\ell u$-group $(G,u)$ is archimedean.

Instead of $\ell u$ groups, one may consider $\ell$-rings, Riesz spaces(vector lattices) or \textit{f}-algebras \cite{Birk,BP} with strong unit and axiomatize the unit interval. The structures obtained in this manner have an MV-algebra reduct endowed with a product operation which can be internal or external.

{\em Product MV-algebras} (PMV-algebras for short) have been defined in \cite{DiND} in the general case and in \cite{MonPMV} in a slightly different way for the unital and commutative case. They are MV-algebras endowed with a binary internal product that satisfies the following, for any $x,y,z \in P$:

(PMV1) $c\cdot (a\odot (a\wedge b)^*)=(c\cdot a)\odot (c\cdot (a\wedge b))^*$

(PMV2) $(a\odot (a\wedge b)^*) \cdot c=(a\cdot c)\odot ((a\wedge b)\cdot c)^*$.

(PMV3) $a \cdot (b\cdot c)= (a \cdot b)\cdot c$.\\
A PMV-algebra is unital if it has a unit for the product, and a PMV\textit{f}-algebra \cite[Theorem 5.4]{DiND} is a PMV-algebra that satisfies the \textit{f}-property:

(f) if $a\wedge b=0$, then $(a\cdot c)\wedge b=(c\cdot a)\wedge b=0$, for any $a,b,c\in P$.\\
It is straightforward that unital PMV-algebras are PMV\textit{f}-algebras, and any PMV\textit{f}-algebra is subdirect product of totally ordered PMV-algebras \cite[Proposition 5.5]{DiND}.

Let us denote by $\mathbf{PMV}$ and $\mathbf{uR}$ the categories of PMV-algebras and $\ell u$-rings such that $u \cdot u \le u$ with suitable morphisms.  In \cite{DiND} the functor $\Gamma$ was extended to a  functor $\Gamma_{(\cdot)}: \mathbf{uR}\rightarrow \mathbf{PMV}$  which is also an equivalence. 

A further extension of the notion of MV-algebra has been introduced in \cite{LeuRMV}. A {\em Riesz MV-algebra}  is a structure $(R, \oplus, ^*, \{\alpha\mid \alpha\in [0,1]\},0)$ such that 
$(R, \oplus,  ^*, 0)$ is an MV-algebra and for any $\alpha,\beta\in [0,1]$ and any $a,b\in R$ we have

(RMV1) $\alpha (a\odot b^*)= (\alpha a)\odot (\alpha b)^*$,

(RMV2) $\max(0,\alpha-\beta) a=(\alpha a)\odot (\beta a)^*$,

(RMV3) $\alpha  (\beta a)= (\alpha\beta) a$,

(RMV4) $1 a=a$.\\
Any  homomorphism of MV-algebras between Riesz MV-algebras preserves the additional unary operations, so it is a homomorphism of Riesz MV-algebras.  Riesz MV-algebras are, up to isomorphism,  unit intervals of Riesz spaces with strong unit. 
Let us denote by $\mathbf{RMV}$ and $\mathbf{uRS}$ the categories of Riesz MV-algebras and, respectively,  Riesz spaces with suitable morphisms.  In \cite{LeuRMV} the functor $\Gamma$ was extended to a  functor $\Gamma_{\mathbb{R}}: \mathbf{uRS}\rightarrow \mathbf{RMV}$  which is also an equivalence.

Finally, {\em \textit{f}MV-algebras} are introduced in \cite{LL1} as algebraic structures $(A, \oplus, \ ^*, \cdot, \{ \alpha \}_{\alpha \in [0,1]} , 0)$ such that $(A, \oplus, \ ^*, \cdot ,0)$ is a PMV\textit{f}-algebra, $(A, \oplus, \ ^*, \{ \alpha \}_{\alpha \in [0,1]} , 0)$ is a Riesz MV-algebra and the condition $\alpha (a\cdot b)=(\alpha b)\cdot b=a\cdot (\alpha b)$ is satisfied for any $\alpha \in [0,1]$ and $a,b\in A$.  If there exist a unit for the product, $A$ will be called unital. The corresponding lattice-ordered structures are the \textit{f}-algebras with strong unit \cite{BP}. If we denot by $\mathbf{fMV}$ and $\mathbf{fuAlg}$ the categories of \textit{f}MV-algebras and \textit{fu}-algebras with suitable morphisms respectively, we establish a categorical equivalence 
$\Gamma_f:\mathbf{fuAlg}\to \mathbf{fMV}$  and, in this case, the functor $\Gamma_f$ extends the previous ones:
$\Gamma$, $\Gamma_{(\cdot)}$, $\Gamma_{\mathbb R}$.

In order to summarize all categorical equivalence, we present the following diagram, in which all horizontal arrows are suitable forgetful functors.

\begin{center}
\begin{tikzpicture}
  \node (A) {$\mathbf{uR}$};
  \node (B) [below of=A] {$\mathbf{PMV}$};
  \node (C) [right of=A] {$\mathbf{auG}$};
  \node (D) [below of=C] {$\mathbf{MV}$};
  \node (E) [right of=C] {$\mathbf{uRS}$};
  \node (F) [below of=E] {$\mathbf{RMV}$};
\node(H)[left of=A] {$\mathbf{fuAlg}$};
\node(G)[below of=H] {$\mathbf{fMV}$};
\node(J)[right of=E] {$\mathbf{fuAlg}$};
\node(I)[below of=J] {$\mathbf{fMV}$};
 \draw[->] (A) to node [swap] {$\Gamma_{(\cdot)}$} (B);
\draw[->] (H) to node [swap] {$\Gamma_{f}$} (G);
   
\draw[->] (J) to node[swap] {$\mathcal{U}_{(\cdot \ell)}$} (E);
 \draw[->] (A) to node {$\mathcal{U}_{(\cdot \ell)}$} (C);
  \draw[->] (C) to node [swap] {$\Gamma$} (D);
  \draw[->] (B) to node [swap] {$\mathcal{U}_{(\cdot)}$} (D);
\draw[->] (I) to node {$\mathcal{U}_{(\cdot)}$} (F);  
\draw[->] (E) to node [swap] {$\mathcal{U}_{( \ell \mathbb{R})}$} (C);
\draw[->] (H) to node {$\mathcal{U}_{( \ell \mathbb{R})}$} (A);  
\draw[->] (E) to node {$\Gamma_{\mathbb{R}}$} (F);
\draw[->] (G) to node[swap] {$\mathcal{U}_{\mathbb{R}}$} (B);
    \draw[->] (J) to node {$\Gamma_{f}$} (I);
  \draw[->] (F) to node {$\mathcal{U}_{\mathbb{R}}$} (D);
 \node [below=0.5cm, align=flush center] at (D){ Figure 1.};
\end{tikzpicture}
\end{center}

We finally mention that a PMV-algebra, a Riesz MV-algebra or a {\em f}MV-algebra is semisimple if its MV-algebra reduct is  semisimple MV-algebra. 

\subsection{States. The  state-completion.}

The notion of \textit{state} for an MV-algebra has be introduced in \cite{Mu}, in relation to the notion of "average degree of truth" of a proposition. See also \cite{MunBook, MunRie} for advanced topics.
\bdfn
A state is a linear function  $s:A \rightarrow [0,1]$ such that $s(1)=1$.
\edfn

A state $s$ is \textit{faithful} if $s(a)=0$ implies $a=0$ for any $a\in A$. We remark that if there exists a faithful state on $A$, than $A$ is semisimple.

A state for an $\ell u$-group $(G,u)$ is a positive normalized additive map $t:G\rightarrow \mathbb{R}$. By \cite{MunBook} any state defined on an MV-algebra $A$ can be uniquely extended to a state on the corresponding $\ell u$-group.

For states in a MV-algebra an equivalent form of the Riesz representation theorem holds, due to Kroupa and Panti \cite{K, P}. 
\bthm \label{teo:KP}
For any MV-algebra $A$ there is an affine isomorphism $v\mapsto s_v$ of the convex set of regular Borel probability measures on the maximal spectral space $Max(A)$ onto the set of states on $A$. For every $f\in A$ and $m\in Max(A)$,
$$ s_v(f)=\int_{Max(A)}f^*(m)dv(m)$$
where $f\mapsto f^*$ is the representation for semisimple MV-algebras by continuous functions.
\ethm

\bfact \label{rem:st01}
The notion of states extends to PMV-algebras, Riesz MV-algebra and \textit{f}MV-algebras without changes in the the definition. In \cite{LeuRMV} the authors prove that for Riesz MV-algebras any state is homogeneous, i.e. it preserves the scalar product.
\efact

A state is {\em $\sigma$-continuous} if $\lim_ns(a_n)=s(a)$ for any $a_1\leq\cdots\leq a_n\leq\cdots$ in  $A$ such that $\bigvee_na_n=a$.  A pair $(A,s)$ with $A$ a $\sigma$-complete MV-algebra and $s$ a faithful $\sigma$-continuous  state is called \textit{probability MV-algebra} \cite{MunRie}.

The metric completion of an MV-algebra with respect to the metric induced by a state was studied in  \cite{LeuStateC}. We remind it in the sequel, since it is important for the present investigation. 

 The starting point is the remark that, given an MV-algebra $A$ and  a state $s:A\to [0,1]$, one can define a pseudo-metric on $A$ by   $\rho_s:A\times A \rightarrow [0,1]$ defined by $\rho_s(x,y)=s(d(x,y))$  for any $x,y\in A$ \cite{MunRie}. The  pseudo-metric $\rho_s$  is a metric iff $s$ is a faithful state.  

We say that $(A^c, s^c)$ is 
{\em the state-completion} of $(A,s)$ if  $A^c$ is the Cauchy completion $A$  w.r.t. the pseudo-metric $\rho_s$ and 
$s^c([\{a_n\}_n])=\lim_n s(a_n)$ for any Cauchy sequence $\{a_n\}_n$ in $A$. We define $\varphi_A: A\rightarrow A^c$ by 
$\varphi_A(a)=[\{a\}]$ for any $a\in A$.

\bthm \cite{LeuStateC} \label{teo:sc}
Let $(A,s)$ be an MV-algebra with a state and $(A^c, s^c)$ be its state-completion then the following hold:

1) $A^c$ is $\sigma$-complete, $s^c$ is a $\sigma$-continuous faithful state and $s^c\circ\varphi_A=s$,

2) $\varphi_A$ is an embedding iff $s$ is faithful,

3) (Universal Property) For any MV-algebra $C$, for any faithful state $m$ such that $C$ is $\rho_m$-complete and for any state-preserving homomorphism of MV-algebras $f:A \rightarrow C$ there exists a unique state-preserving embedding of MV-algebras $f^c:A^c \rightarrow C$ such that $f^c \circ \varphi_A =f$.

\ethm
\bfact
We remark that any $\sigma$-complete MV-algebra is semisimple \cite[Proposition 6.6.2]{Mu}, and  $(A^c, s^c)$ is a probability MV-algebra. 
\efact

\bprop \label{pro:scPMV}
If $P$ is a unital and commutative PMV-algebra and $s$ a state, $P^c$ is a PMV-algebra.
\eprop
\begin{proof}
In order to define the product on $P^c$, it is enough to prove that the internal product on $P$ is continuous with respect to $\rho_s$. Following the definition in \cite{LeuStateC}, we define 
$[\mathbf{x}]\cdot [\mathbf{y}]=[(x_n\cdot y_n)_n]$ whenever  $\mathbf{x}=\{x_n\}_n$ and $\mathbf{y}=\{y_n\}_n$.
The product is well defined if and only if $\mathbf{x}\sim \mathbf{y}$ implies $\mathbf{x}\cdot \mathbf{z}\sim \mathbf{y}\cdot \mathbf{z}$.
By definition this holds if and only if $\rho_s(x_n\cdot z_n, y_n \cdot z_n)\rightarrow 0$. By property of the unitary product, see \cite[Corollary 5.7]{FloLeu} 
\vsp{-0.2}
\begin{center}
$\rho_s(x_n\cdot z_n, y_n \cdot z_n)= s(d(x_n\cdot z_n, y_n \cdot z_n))=s(z_n \cdot d(x_n, y_n))\le s(d(x_n, y_n))=\rho_s(x_n, y_n) \rightarrow 0$ by hypothesis,
\end{center}
\vsp{-0.2}
therefore the conclusion follows.
\end{proof}

Following \cite{LeuSCRMV}, a \textit{state-complete} Riesz MV-algebra is a structure $(A,s)$ such that $A$ is a Riesz MV-algebra, $s$ is a state on $A$ and $(A, \rho_s)$ is a complete metric space. An $L$-space \cite{JVanR} is a Banach lattice $(L, \| \cdot \|)$ such that 
\vsp{-0.2}
\begin{center}
$x,y\ge 0$ in $L$ implies $\| x+y\|=\|x\|+\|y\|$.
\end{center}
\vsp{-0.2} 

By  \cite[Corollary 1]{LeuSCRMV}  any state-complete Riesz MV-algebra is  Dedekind-MacNeille complete.
In \cite{LeuSCRMV} the author proves the categorical duality between state-complete Riesz MV-algebras and a particular class of measure space. We main point is the following theorem, which is similar with Kakutani representation for abstract $L$-spaces  \cite{Kak}

\bdfn If $(X, \Omega, \mu)$ is a measure space then we define
\vsp{-0.2}
\begin{center}
$L_1(\mu)_u=\{ f\in L_1(\mu) \mid \mathbf{0}\le f\le \mathbf{1}\}$.  
\end{center}
\vsp{-0.2} 
\edfn

\bfact\label{rem:fmv}
We note  that $L_1(\mu)$ is  an \textit{f}-algebra and $\mathbf{1}$ is a weak unit  of  $L_1(\mu)$.  Therefore $L_1(\mu)_u$ is an \textit{f}MV-algebra.
\efact

\bthm \cite{LeuSCRMV} \label{teo:equivL1}
For any state complete Riesz MV-algebra $(A,s)$ there exists a measure space $(X, \Omega, \mu)$  and a  isomorphism of Riesz MV-algebras $I_A:A\to L_1(\mu)_u$  such that $s(a)=\int I_A(a)d\mu$ for any $a\in A$.  Moreover,  $(X, \Omega, \mu)$ is a probability space  such that  $X$ is an extremally disconnected compact Hausdorff space,
$\Omega$ is the Borel $\sigma$-algebra of $X$ and $\mu$ is a topological finite measure.
\ethm

\bdfn \cite{LeuSCRMV}
 A measure space  $(X, \Omega, \mu)$  that satisfies the properties from Theorem \ref{teo:equivL1} is called {\em $L$-space}.
\edfn

\section{Embedding in $L^1(\mu)_u$} \label{section:emb}

We  prove that any MV-algebra which has a faithful state can be embedded in an \textit{f}MV-algebra  of integrable functions. 
As a preliminary step, we embed the MV-algebra in its divisible hull.

\bfact \label{rem:div}
\begin{enumerate}
\item Any MV-algebra can be embedded in a divisible one. For details on divisible MV-algebras see \cite{Brunella}. In the semisimple case, $A^d=\{ a\in C(X) \mid a=\frac{a_1}{n}+ \ldots + \frac{a_n}{n},\ a_i \in A, n\in \mathbb{N} \}$, where $A\subseteq C(X)$. Moreover,  if $P$ is a unital and commutative PMV-algebra then  $P \hookrightarrow P^d$ is an embedding of PMV-algebras.
\item If $s:A\rightarrow [0,1]$ is a state on $A$ then  $s$  can be extended to a state $s^d: A^d \rightarrow [0,1]$ \cite[Theorem 6]{KExt} such that  $s^d(\alpha a)=\alpha s^d(a)$ for any $\alpha \in [0,1] \cap \mathbb{Q}$ and $a\in A$ 
\cite[Lemma 11]{LeuRMV}. Note that   $s^d$ is  faithful whenever $s$ is faithful.

\end{enumerate}
\efact

\bthm \label{teo:embMeasure}
 Let $A$ be an MV-algebra and   $s:A \rightarrow [0,1]$ a state on A. There exists an  $L$-space $(X, \Omega, \mu)$  and  a homomorphism of MV-algebras $F_A: A\rightarrow L_1(\mu)_u$ such that $s(a)=\int F_A(a)d\mu$ for any $a\in A$.  If $s$ is faithful then  $F_A$ is an embedding.
\ethm

\begin{proof} 
 Let $B= \quot{A}{Rad(A)}$, where $Rad(A)$ is the intersection of all the maximal ideals of $A$. Denote by   $\pi: A \rightarrow \quot{A}{Rad(A)}$ the canonical epimorphism. We notice that $B$ is semisimple and we can define the divisible hull $B^d$ as in Remark \ref{rem:div}. Let $\iota_d$ be the embedding in the divisible hull. With the notation of Theorem \ref{teo:sc}, we have
\vsp{-0.2}
\begin{center}
\begin{tikzpicture}
 \node (A) {$A$};
  \node (B) [right=1cm of A] {$B$};
  \node (C) [right=1cm of B] {$B^d$};
  \node (D) [right of=C] {$B^{dc}$};
  \draw[->] (A) to node {$\pi$} (B);
  \draw[right hook->] (B) to node {$\iota_d$} (C);
  \draw[right hook->] (C) to node {$\varphi_{B^{d}}$} (D);
\end{tikzpicture}
\end{center}
\vsp{-0.2}
By \cite[Lemma 3.1]{DiaLeu}, $B^{dc}$ is a state-complete Riesz MV-algebra, therefore by Theorem \ref{teo:equivL1}, $B^{dc}\simeq L_1(\mu)_u$ for a suitable $L$-space.
Finally, $t=\pi \circ s$ is a state on $B$. By Remark \ref{rem:div} $t$ extends to $t^d$, and by Theorem \ref{teo:sc} it extends to $t^{dc}$. The conclusion follows from Theorem \ref{teo:equivL1}.
\end{proof}

\bfact The  integral representation of  a state from Theorem \ref{teo:embMeasure} is obviously  Kroupa-Panti's result \cite{K,P}. We mention that, for our development, the representation of the algebraic structures is also crucial.  We also mention that we followed the approach from \cite{JVanR}, where Riesz integral representation  is derived as a consequence of Kakutani's representation for $L$-spaces.
\efact

\bprop\label{prop:embMeasure}
If $P$ is a unital and semisimple PMV-algebra  (\textit{f}MV-algebra) then the morphism $F_P$ from Theorem \label{teo:embMeasure} is a morphism of PMV-algebras (\textit{f}MV-algebras). 
\eprop
\begin{proof}
 We first remark that any unital and semisimple PMV-algebra or \textit{f}MV-algebra is commutative by \cite{LL1}.
If $A$ is a unital and semisimple PMV-algebra, the conclusion follows by Proposition \ref{pro:scPMV}, Remark \ref{rem:div} and  \cite{Conrad}, since in any  archimedean \textit{f}-rings, the ring structure is generated by the additive group, and therefore any  homomorphism of groups for the group reduct is an homomorphism of rings, and the same applies for unital and semisimple PMV-algebras.  If $A$  is an unital and semisimple \textit{f}MV-algebra, the result  follows by Remark \ref{rem:fmv}, \cite[Proposition 3.2]{LL1} and \cite[Corollary 2]{LeuRMV}. In particular in \cite{LL1} is proved that any linear homomorphism between unital and semisimple \textit{f}MV-algebras commutes with the internal product.
\end{proof} 

By means of this representation, we prove two  immediate consequences for states  using  results in functional analysis:
H\"{o}lder's inequality and the Hausdorff moment problem.

\medskip

\noindent{\bf  H\"{o}lder's inequality for PMV-algebras and \textit{f}MV-algebras}

\noindent The first result  H\"{o}lder's inequality for PMV-algebras and \textit{f}MV-algebra, in the unital and semisimple case. We recall that any unital and semisimple algebra is commutative.
\bthm
Let $A$ be a semisimple $PMV^+$-algebra ($\mathbb{FR}^+$-algebra) and  $s:A\to [0,1]$ a state. If $p$, $q\in [1,\infty)$ with $\frac{1}{p}+\frac{1}{q}=1$ then 

\begin{center}
$s(a\cdot b)\leq s(a^p)^{\frac{1}{p}} s(b^q)^{\frac{1}{q}}$ for any $a$, $b\in A$.
\end{center}
\vsp{-0.2}
\ethm
\begin{proof}
By Proposition \ref{prop:embMeasure}, $F_A:A\to L_1(\mu)_u$  is a morphism of  PMV-algebras (\textit{f}MV-algebras). By H\"{o}lder's inequality for $L_1(\mu)$, for any  $a$, $b\in A$ we get

\begin{center}
$\int_X F_A(a\cdot b)d\mu=\int_X(F_A(a)\cdot F_B(b))d\mu\leq \left( \int_X F_A(a^p)d\mu\right) ^{\frac{1}{p}} \left( \int_X F_A(b^q)d\mu\right) ^{\frac{1}{q}},$
\end{center}

\noindent and by Theorem \ref{teo:equivL1},  $s(a\cdot b)\leq s(a^p)^{\frac{1}{p}} s(b^q)^{\frac{1}{q}}$.
\end{proof}
\quad

\medskip 
\noindent{\bf   The Hausdorff moment problem for PMV-algebras and  \textit{f}MV-algebras}

\noindent In statistics and probability a very central subject is the Moment Problem. Given a interval $I\subseteq \mathbb R$, the $n^{th}$-moment of a probability measure $\mu$ on $I$ is defined as $\int _{I}x^n d\mu $. Let $\{ m_k \}_{k\ge 0}$ be a sequence of real numbers, the Moment Problems on $I$ consists on finding out the condition on $\{ m_k \}_{k\ge 0}$ for which there exists a probability measure $\mu$ on $I$ such that $m_k$ is the $k^{th}$ moment of $\mu$.\\
When $I=[0,1]$ we get the Hausdorff moment problem \cite{Ha1,Ha2}. We will prove a  similar result in the context of MV-algebras. 

For any $k\geq 1$ we define  $p_k:[0,1]\rightarrow [0,1]$  by $p_k(x)=x^k$ for any $x\in [0,1]$. We also set 
$p_0(x)=1$ for any $x\in [0,1]$. Note that $p_k\in FR_1$ for any $k\geq 0$.\\

\noindent If  $\{m_k | k\geq 0\} $  a sequence of real numbers in $[0,1]$, we define:
\vsp{-0.2}
\begin{center}
$\Delta ^0 m_k = m_k$, \quad $\Delta ^r m_k = \Delta ^{r-1} m_{k+1} -\Delta^{r-1} m_k $  for any $r$, $k\geq 0$.
\end{center}
\vsp{-0.2}
The sequence $\{ m_k \}_k$ satisfies the {\em Hausdorff moment condition} if $m_0=1$  and $(-1)^r\Delta ^r m_k \ge 0$ for any $r$, $k\geq 0$
\cite{Fe}.

\bthm \label{teo:01MP}
Let $C$ be any unital and semisimple PMV-subalgebra (unital and semisimple \textit{f}MV-subalgebra) of $C([0,1])$ such that $p_1\in C$. There exists a state $s:C\rightarrow [0,1]$ such that $s(p_k)=m_k$ if and only if the sequence $\{ m_k \}$ satisfies the {\em Hausdorff moment condition}.
\ethm
\begin{proof}
Let $s$ be a state such that $s(x^k)=m_k$. Since $C$ is unital, the set of its ideals coincide with the set of ideals of its MV-algebra reduct and by general theory of MV-algebras (see for example \cite{CDM}, Chapter 3) and by the integral representation due to Kroupa and Panti \cite{K, P} we have
$$s(f)= \int_{0}^1 f d\mu,$$
for any $f\in C$, where $\mu$ is a probability measure on $[0,1]$.\\
By \cite{Fe} we have $$(-1)^r\Delta ^r m_k= \sum_{h=0}^r \binom{r}{h}(-1)^h m_{k+h},$$
then by the hypothesis,
$$(-1)^r \Delta ^r m_k= \sum_{h=0}^r (-1)^h \binom{r}{h}\int_{0}^1 x^{k+h}d\mu = \int_{0}^1 \left[ x^k \sum_{h=0}^r \binom{r}{h} (-1)^h x^h \right]  d\mu =$$ $$= \int_{0}^1 x^k (1-x)^r d\mu \ge 0,$$
therefore the Hausdorff moment condition is satisfied.\\
On the other hand, let $\overline s$ the functional on the set $\{ p_n \mid n=0,1,2,\ldots \}$ such that $\overline s (p_k)=m_k$. By \cite{Mir} $\overline s$ has a unique extension $\widetilde s$ to a linear prevision (that is a positive and normalized linear functional) from $C([0,1], \mathbb R)$ to $\mathbb R$. In particular, $ \widetilde s $ is a state between $\ell$-groups, then $s: \Gamma ( C([0,1], \mathbb R), 1) \rightarrow [0,1]$ is a state on $C([0,1])$ (see for example \cite{MunBook}). Taking $s|_{C}$, the restriction of $s$ to $C$ we get the desired result. 
\end{proof}

\noindent Let $A$ be a unital PMV-algebra (unital \textit{f}MV-algebra) such that $\quot{A}{Rad(A)} \subseteq C([0,1])$, and let $\phi _A$ be the map $\phi _A : A \to C([0,1])$ obtained composing the canonical epimorphism $A\to \quot{A}{Rad(A)}$ and the embedding of $\quot{A}{Rad(A)}$ in $C([0,1])$. Moreover, we ask that $p_1 \in \phi _A (A)$.
\bcor
Let A be a unital PMV-algebra (unital \textit{f}MV-algebra) as defined above. If the sequence $\{ m_k \}$ satisfies the {\em Hausdorff moment condition}, then there exists a state $s:A\rightarrow [0,1]$ such that $s(p_k)=m_k$.
\ecor
\begin{proof}
Theorem  \ref{teo:01MP} holds for $\quot{A}{Rad(A)}$, therefore $s$ will be the composition of the state $t: A/Rad(A) \to [0,1]$ with the canonical epimorphism.
\end{proof}
\quad

\section{Stochastic independence of probability MV-algebras} \label{section:ind}

The notion of independence for probability MV-algebra is one of the open problems mentioned in \cite{MunRie}. A partial solution has given in \cite{LeuStateC} for MV-algebras with states, MV-algebras with extremal states and semisimple MV-algebras. In this section we obtain a solution for probability MV-algebras.

Recall that a probability MV-algebra is  a pair $(A,s)$ with $A$ a $\sigma$-complete MV-algebra and $s$ a faithful $\sigma$-continuous  state.

\bdfn
Let $(A, s_A)$, $(B, s_B)$ and $(T,s)$ probability MV-algebras, and $\beta: A\times B\rightarrow T$ a bilinear function. $(A, s_A)$ and $(B, s_B)$ are said to be $(T, s, \beta)$-\textit{independent} if 
$s(\lambda(a,b))=s_A(a)\cdot s_B(b)$, for any $a\in A$ and $b\in B$.
\medskip

Given   $(A, s_A)$ and $(B, s_B)$  two probability MV-algebras, our problem is to define a probability MV-algebra $(T,s)$  and a bilinear  function $\beta: A\times B\rightarrow T$ such that  $(A, s_A)$ and $(B, s_B)$ are  $(T, s, \beta)$-independent.
\edfn

\bfact\label{pr}
We recall some general results from Measure Theory, see \cite[253D, 253G, 253F]{Fre} for further details.
Let $(X_A,\Omega_A,\mu_A)$ and  $(X_B,\Omega_B,\mu_B)$ be measure spaces. There is a measure space $(X_A\times X_B,\Lambda,\lambda)$   where $\lambda$ is the c.l.d. product measure on $X_A\times X_B$ and the following properties hold:\\
(1)  $\otimes:L_1(\mu_A)\times L_1(\mu_B)\to L_1(\lambda)$, $(f,g) \mapsto f\otimes g$  is a bounded bilinear operator,\\
(2)  $\int{(f\otimes g)}d\lambda=\int{f}d\mu_A\int{g}d\mu_B$ whenever $f\in L_1(\mu_A)$, $g\in L_1(\mu_B)$,\\
(3) $f\otimes g\geq 0$ in $L_1(\lambda)$ whenever $f\geq 0$ and $g\geq 0$,\\
(4) the following universal property is satisfied: for any  Banach lattice $W$ (norm complete Riesz space)  and  bilinear function $\phi$  there exists a unique linear function $\omega$ such that $\omega\circ\otimes=\phi$.
\begin{center}
\begin{tikzpicture}
 \node (A) {$L_1(\mu_A) \times L_1(\mu_B)$};
  \node (B) [right=1cm of A] {$L_1(\lambda )$};
  \node (C) [below of=A] {$W$};
  \draw[->] (A) to node {$\otimes$} (B);
  \draw[->] (A) to node [swap] {$\phi$} (C);
  \draw[dashed, ->] (B) to node {$\omega$} (C);
\end{tikzpicture}
\end{center}
\vsp{-0.2}
\efact

\bdfn\label{defT}
For any $(A, s_A)$ and $(B, s_B)$ probability MV-algebras and let $(X_A \Omega_A, \mu_A)$, $(X_A \Omega_A, \mu_A)$ be defined by Theorem \ref{teo:equivL1} and $(X_A \times X_B , \Lambda, \lambda)$ is the product measure space from Remark \ref{pr}. The  \textit{product space} is $(T, s_T)$ where
\begin{center}
 $T=L_1(\lambda)_u$ and ${s_T(f)=\int_{X_A\times X_B} f d\lambda}$ for any $f\in T$.
\end{center}

Assume $\beta:A \times B \rightarrow L_1(\lambda)_u$ is the bilinear map defined by 
\vsp{-0.2}
\begin{center}
$\beta(a,b)=f_a\otimes f_b$, where

$A\hookrightarrow A^d \hookrightarrow A^{dc} \simeq L_1(\mu_A)_u$, $a\mapsto F_A(a)= f_a$, 

$B\hookrightarrow B^d \hookrightarrow B^{dc} \simeq L_1(\mu_B)_u$, $b\mapsto F_B(b)=f_b$.
\end{center}
\vsp{-0.2}
\edfn

\bfact
With the above definition, $(T, s_T)$ is a probability MV-algebra by Theorem \ref{teo:equivL1} and Theorem \ref{teo:sc}. Moreover,  by Theorem \ref{teo:equivL1} and \cite[253D]{Fre},
\vsp{-0.2}
\begin{center}
$s_T(\beta (a,b))=\int (f_a\ten f_b d\lambda)=\int f_ad\mu_A \int f_b d\mu_B$.
\end{center}
\vsp{-0.2}
\efact

\bdfn
Let $(A, s_A)$, $(B, s_B)$ and $(C, s_C)$  be probability MV-algebras. A linear function $\omega:A\to C$ is {\em bounded} if 
 there exists a natural number  $K\geq 1$ such that 

\vsp{-0.2}
\begin{center}
$s_C(\omega(a))\leq K_\oplus s_A(a)$ for any $a\in A$. 
\end{center}
\vsp{-0.2}

We say that a bilinear function  $\gamma:A \times B \rightarrow C$ is {\em bounded} if  there exists a natural number  $K\geq 1$ such that 

\vsp{-0.2}
\begin{center}
$s_C(\gamma(a,b))\leq K_\oplus (s_A(a)s_B(b))$ for any $a\in A$ and $b\in B$. 
\end{center}
\vsp{-0.2}

The bilinear function  $\gamma$ is {\em continuous} if the following property holds: for any 
$(x_n)_n\subseteq A$, $x\in A$, $(y_n)_n\subseteq B$, $y\in B$, 
\vsp{-0.2}
\begin{center}
$\rho_{s_A}(x_n,x) \rightarrow 0$ and $\rho_{s_B}(y_n,y) \rightarrow 0$ imply  $\rho_{s_C}(\gamma(x,y),\gamma (x_n, y_n)) \rightarrow 0$.
\end{center}
\vsp{-0.2}
\edfn  

One can immediately see that any bounded linear function is continuous. In the sequel we prove the same result for bilinear functions.

\blem\label{lem:cont}
If $(A, s_A)$, $(B, s_B)$ and $(C, s_C)$  are  probability MV-algebras,  $K\geq 1$  is a  natural number and   the bilinear function  $\gamma :A \times B \rightarrow C$ is $K$- bounded then, for any $a$, $a^\prime\in A$ and $b$, $b^\prime\in B$
\vsp{-0.2}
\begin{center}
$\rho_{s_C}(\gamma(a, b), \gamma^d(a^\prime, b^\prime))\leq K_\oplus (\rho_{s_A}(a, a^\prime)\oplus  \rho_{s_B}(b, b^\prime))$.
\end{center}
\vsp{-0.2}
\elem
\begin{proof}
Assume $b\in B$.  Since $\gamma(\cdot,b):A\to C$ is a linear map, we get
\vsp{-0.2}
\begin{center}
$d(\gamma (a,b),\gamma (a^\prime,b))\leq \gamma (d(a, a^\prime),b)$ for any $a$, $a^\prime\in A$.
\end{center}
\vsp{-0.2}
It follows that  
\vsp{-0.2}
\begin{center}
$s_C(d(\gamma(a,b),\gamma(a^\prime,b))\leq s_C( \gamma(d(a, a^\prime),b))\leq K_\oplus ( s_A(d(a, a^\prime)) s_B(b))$,
\end{center}
\vsp{-0.2}
so $s_C(d(\gamma(a,b),\gamma(a^\prime,b))\leq K_\oplus s_A(d(a, a^\prime))$ 
 for some constant $K\geq 0$.

 Similarly we get $s_C(d(\gamma(a,b),\gamma(a,b^\prime))\leq K_\oplus s_B(d(b, b^\prime))$ for any $a\in A$, so

\vsp{-0.2}
\begin{center}
  $s_C(d(\gamma(a,b),\gamma(a^\prime,b^\prime)))\leq
       s_C(d(\gamma(a,b),\gamma(a^\prime,b))) \oplus s_C(d(\gamma(a^\prime,b),\gamma(a^\prime,b^\prime)))$\\
$ \leq  K_\oplus  s_A(d(a, a^\prime)) \oplus K_\oplus s_B(d(b, b^\prime))$, so\\
$\rho_{s_C}(\gamma(a, b), \gamma^d(a^\prime, b^\prime))\leq K_\oplus (\rho_{s_A}(a, a^\prime)\oplus  \rho_{s_B}(b, b^\prime))$.

\end{center}
\vsp{-0.2}
\end{proof}

\bcor\label{cor:cont}
If $(A, s_A)$, $(B, s_B)$ and $(C, s_C)$  are  probability MV-algebras then any bounded bilinear  function  $\gamma :A \times B \rightarrow C$ is continuous. 
\ecor
\begin{proof}
It follows by Lemma \ref{lem:cont}.
\end{proof}

\bprop\label{pr:extlinfct}
If  $A$ and $B$ are  MV-algebras and $\sigma :A \rightarrow B$ is a linear function, then  there is a unique linear function $\sigma^d: A^d \rightarrow B^d$ that extends $\sigma$.
\eprop
\begin{proof}
Let $a\in A^d$ and  $a_1,\ldots, a_n\in A$ such that $a=\frac{a_1}{n}+\cdots +\frac{a_n}{n}$. We set
$\sigma^d(a)=\frac{1}{n}(\sigma(a_1)+\cdots +\sigma(a_n))$. Using Riesz decomposition property\cite[Section 2.9]{LDHandbook}  in $A^d$ one can prove that $\sigma^d$ is well-defined. We show that $\sigma^d$ is linear. Assume that $a+a^\prime$ is defined in $A^d$. We know that  $a=\frac{a_1}{n}+\cdots +\frac{a_n}{n}$ and  $a^\prime=\frac{a^\prime_1}{m}+\cdots +\frac{a^\prime_m}{m}$ where 
 $a_1,\ldots, a_n, a_1^\prime,\ldots, a_n^\prime\in A$. It follows that 
\vsp{-0.2}
\begin{center}
$a+a^\prime= \frac{a_1}{n}+\cdots +\frac{a_n}{n} +\frac{a^\prime_1}{m}+\cdots +\frac{a^\prime_m}{m}=m\frac{a_1}{nm}+\cdots +m\frac{a_n}{nm} +n\frac{a^\prime_1}{nm}+\cdots +n\frac{a^\prime_m}{nm}$.
\end{center}
\vsp{-0.2}
We get
\vsp{-0.2}
\begin{center}
$\sigma^d(a+a^\prime)=\frac{1}{nm}(m\sigma(a_1)+\cdots +m\sigma(a_n)+n\sigma(a^\prime_1)+\cdots +n\sigma(a^\prime_m))$

$= \frac{1}{n}(\sigma(a_1)+\cdots +\sigma(a_n)) + \frac{1}{m}(\sigma(a^\prime_1)+\cdots +\sigma(a^\prime_m))$

$=\sigma(a)+\sigma(a^\prime).$
\end{center}
\vsp{-0.2}
\end{proof}

\blem \label{lem:extBimDiv}
If  $(A, s_A)$, $(B, s_B)$ and $(C, s_C)$  are probability MV-algebras and $\gamma :A \times B \rightarrow C$ is a bilinear function, then  there exists a unique bilinear $\gamma^d: A^d \times B^d \rightarrow C^d$ that extends $\gamma$. Moreover, if $\gamma$ is bounded then $\gamma^d$ is also bounded. 
\elem
\begin{proof}
 We first recall that any linear function between divisible MV-algebra is linear w.r.t. scalars in $[0,1]\cap \mathbb{Q}$, has remarked for states in Remark \ref{rem:div}.  If $a\in A^d$ and $b\in B^d $ then there are $a_1,\ldots, a_n\in A$ and $b_1,\ldots, b_m\in B$ such that $a=\frac{a_1}{n}+\ldots +\frac{a_n}{n}$ and $b= \frac{b_1}{m}+\ldots +\frac{b_m}{m}$. 
We define $\gamma^d: A^d \times B^d \rightarrow C^d$ by
\vsp{-0.2}
\begin{center}
$\gamma^d (a, b)= \frac{1}{nm}\sum \gamma (a_i, b_j)$.
\end{center}
\vsp{-0.2}
The fact that $\gamma^d$ is well-defined and the uniqueness follow by Proposition \ref{pr:extlinfct}. We have that 
\vsp{-0.2}
\begin{center}
$s_C^d(\gamma^d(a,b))= s_C^d(\gamma^d (\frac{a_1}{n}+\ldots +\frac{a_n}{n}, \frac{b_1}{m}+\ldots +\frac{b_m}{m})) = s_C^d (\frac{1}{nm}\sum_{ij} \gamma (a_i, b_j))= s_C^d (\sum_{ij} \frac{1}{nm}\gamma (a_i, b_j))=\sum_{ij} \frac{1}{nm} s_C^d (\gamma (a_i, b_j)) \le \sum_{ij}\frac{1}{nm} ( K_\oplus (s_A(a_i) s_B(b_j))) = $\\

$ \sum_{i,j} \frac{1}{nm} \min (K s_A(a_i)s_B(b_j), 1) =  \sum_{i,j} \min (K s_A^d(\frac{a_i}{n})s_B^b(\frac{b_j}{m}), \frac{1}{nm}) \le \min (K (\sum_{i,j} s_A^d(a_i)s_B^d(b_j) ), \frac{nm}{nm})=
K_\oplus\left(\sum_{i,j} s_A^d\left( \frac{a_i}{n}\right) s_B^d\left( \frac{b_j}{m}\right) \right)=
K\oplus\left( \left( s_A^d\left( \frac{a_1}{n}\right) +\ldots +s_A^d\left( \frac{a_n}{n}\right) \right) \left( s_B^d\left( \frac{b_1}{m}\right) +\ldots +s_B^d\left( \frac{b_m}{m}\right) \right)\right) = K_\oplus( s_A^d(a) s_B^d (b))$
\end{center}
\vsp{-0.2}
so $\gamma^d$ is bounded.
\end{proof}

\bprop \label{pro:extBimDiv}
Let $(A, s_A)$, $(B, s_B)$ and $(C, s_C)$ be probability MV-algebras, and $\gamma :A \times B \rightarrow C$ a bounded bilinear function. Then   there exists a unique bounded bilinear  function  $\gamma^{c}:A^{c}\times B^{c}\to C^{c}$ that extends $\gamma$, defined by
\vsp{-0.2}
\begin{center}
 $\gamma^{c}([\{a_n\}_n],[\{b_n\}_n])=[\{\gamma(a_n,b_n)\}_n]$ 
\end{center}
\vsp{-0.2}
for any Cauchy sequences  $\{a_n\}_n$ from $A$ and $\{b_n\}_n$ from $B$. 
\eprop

\begin{proof} Let $K\geq 1$ be a natural number such that  $\gamma$ is bounded with constant $K$.
If $\{a_n\}_n$, $\{a_n^\prime\}_n$  are  Cauchy sequences in $A$ and let $\{b_n\}_n$,  $\{b_n^\prime\}_n$ are Cauchy sequences in $B$,   from Lemma \ref{lem:cont},  we infer that $\{\gamma(a_n,b_n)\}_n$ is a Cauchy sequence in $C$. Moreover,
$[\{a_n\}_n]=[\{a_n^\prime\}_n]$ and $[\{b_n\}_n]=[\{b_n^\prime\}_n]$ imply $[\{\gamma(a_n,b_n)\}_n]=[\{\gamma(a_n^\prime,b_n^\prime)\}_n]$.  Hence $\gamma^c$ can be defined by 
\vsp{-0.2}
\begin{center}
 $\gamma^{c}([\{a_n\}_n],[\{b_n\}_n])=[\{\gamma(a_n,b_n)\}_n]$ 
\end{center}
\vsp{-0.2}
for any Cauchy sequences  $\{a_n\}_n$ from $A$ and $\{b_n\}_n$ from $B$.  We get
\vsp{-0.2}
\begin{center}
$s_C^c(\gamma^{c}([\{a_n\}_n],[\{b_n\}_n]))=\lim_n s_C(\gamma(a_n,b_n))\leq K_\oplus (\lim_n s_A(a_n)\lim_n s_B(b_n))$, \\
$s_C^c(\gamma^{c}([\{a_n\}_n],[\{b_n\}_n]))\leq K_\oplus (s_A^c([\{a_n\}_n]) s_B^c([\{b_n\}_n]))$. 
\end{center}
\vsp{-0.2}
Let $\beta:A^{c}\times B^{c}\to C^{c}$ be another bilinear bounded function that extends $\gamma$. Since $\lim_n[a_n]=[\{a_n\}_n]$ in $A^c$ and  $\lim_n[b_n]=[\{b_n\}_n]$ in $B^c$, by Corollary \ref{cor:cont}, 
\vsp{-0.2}
\begin{center}
$\beta([\{a_n\}_n],[\{b_n\}_n])=\lim_n \beta([a_n],[b_n])=\lim_n\gamma(a_n,b_n)=[\{\gamma(a_n,b_n)\}_n]= \gamma^{c}([\{a_n\}_n],[\{b_n\}_n])$.
\end{center}
\vsp{-0.2}
\end{proof}

\bfact
If  $(A, s_A)$, $(B, s_B)$ and $(C, s_C)$  are probability MV-algebras and $\gamma :A \times B \rightarrow C$ is a bounded bilinear function then the following diagram is commutative:
\begin{center}
\begin{tikzpicture}
  \node (A) {$A\times B$};
  \node (B) [right=2.5cm of A] {$A^d\times B^d$};
  \node (D) [right=2.5cm of B] {$A^{dc}\times B^{dc}$};
  \node (C) [below of=A] {$C$};
  \node (E) [below of=B] {$C^d$};
  \node (F) [below of=D] {$C^{dc}$};
  \draw[right hook->] (A) to node {\small{$\iota_{d_A} \times \iota_{d_B}$}} (B);
  \draw[right hook->] (B) to node {\small{$\varphi_{A^d} \times \varphi_{B^d}$}} (D);
  \draw[->] (A) to node {$\gamma$} (C);
  \draw[right hook->] (C) to node {$\iota_{dC}$} (E);
  \draw[right hook->] (E) to node {$\varphi_C$} (F);
  \draw[dashed, ->] (B) to node {$\gamma^d$} (E);
  \draw[dashed, ->] (D) to node {$\gamma^{dc}$} (F);
\end{tikzpicture}
\end{center}
This is a straightforward consequence of Lemma \ref{lem:extBimDiv} and Proposition \ref{pro:extBimDiv}.
\efact

We are ready to prove the main result of this section.

\bthm \label{teo:unpropInd}
 Let $(A, s_A)$, $(B, s_B)$ be probability MV-algebras and assume  $(T,s_T)$ and $\beta:A\times B\rightarrow T$ are defined as in Definition \ref{defT}.
For any probability MV-algebra $(C, s_C)$  and any bounded bilinear  function  $\gamma : A\times B \rightarrow C$ there exists a unique bounded  linear function $\omega: T \rightarrow C^{dc}$ such that $\omega(\beta(a,b))=\gamma(a,b)$ for any $a\in A$, $b\in B$, i.e. 
 $\omega(f_a\ten f_b)=f_{\gamma(a,b)}$ for any $a\in A$, $b\in B$. 

 \begin{center}
\begin{tikzpicture}
  \node (A) {$A\times B$};
  \node (B) [right=1cm of A] {$T=L_1(\lambda)_u$};
  \node (F) [below=2.5cm of A] {$C^{dc}=L_1(\mu_C)_u$};
    \node (C) at ($(A)!0.5!(F)$) {$C$};
  \draw[->] (A) to node {$\beta$} (B);
  \draw[->] (A) to node [swap]{$\gamma$} (C);
  \draw[right hook->] (C) to node [swap]{$\varphi_{C^d}$} (F);
  \draw[dashed, ->] (B) to node {$\omega$} (F);
\end{tikzpicture}
\end{center}
 
\ethm
\begin{proof}

By Proposition \ref{pro:extBimDiv} there exists a unique bounded  bilinear function f$\gamma^{dc}: A^{dc}\times B^{dc} \rightarrow C^{dc}$ that extends $\gamma$. By Theorem \ref{teo:equivL1} $A^{dc}\simeq L_1(\mu_A)_u$, $B^{dc}\simeq L_1(\mu_B)_u$ and $C^{dc}\simeq L_1(\mu_C)_u$ for suitable measure spaces. 

Note that $\mathbf{1}_A$ and  $\mathbf{1}_B$ are weak units in $ L_1(\mu_A)$ and  $L_1(\mu_B)$. If we set
\vsp{-0.2}
\begin{center}
$L_A=\{ f\in  L_1(\mu_A)\mid |f|\leq n {\mathbf 1}_A  \mbox{ for some } n\geq 1 \}$,
$L_B=\{ g\in  L_1(\mu_B)\mid |g|\leq n {\mathbf 1}_B  \mbox{ for some } n\geq 1 \}$,
\end{center}
\vsp{-0.2}
then $(L_A,{\mathbf 1}_A)$  and $(L_B,{\mathbf 1}_B)$ are \textit{f}-algebras with strong unit and they are dense in $ L_1(\mu_A)_u$ and  $L_1(\mu_B)_u$, respectively.  By \cite[Proposition 6.5]{FloLeu} there exists an extension  $\widetilde{\gamma}: L_A \times L_B \rightarrow L_1(\mu_C)$ \textit{u}-bilinear function that extends $\gamma^{dc}$.  We recall that by construction $\widetilde{\gamma}(\mathbf{1}_A, \mathbf{1}_B)\le \mathbf{1}_C$. One can easily see that the bilinear map $\widetilde{\gamma}$  is bounded.

In order to extend $\widetilde{\gamma}$  to $L_1(\mu_A)\times L_1(\mu_B)$ we shall apply the B.L.T. theorem \cite[Theorem I.7]{BLT}
twice.   Let $g\in L_B$ be an arbitrary element and apply the B.L.T. theorem for $\widetilde{\gamma}(\cdot, g):L_A\to  L_1(\mu_C)$. Hence there exists a unique bounded linear transformation $\gamma_g:  L_1(\mu_A)\to   L_1(\mu_C)$ such that 
$\gamma_g(f)=\widetilde{\gamma}(f,g)$ for any $f\in  L_1(\mu_A)$. Now, we fix  $f\in  L_1(\mu_A)$ and we define
$\gamma^f:L_B\to  L_1(\mu_C)$ by $\gamma^f(g)=\gamma_g(f)$ for any $g\in L_B$. Applying again the B.L.T. theorem  we get a unique  bounded linear transformation ${\gamma^f}^\prime:  L_1(\mu_B)\to  L_1(\mu_C)$ such that ${\gamma^f}^\prime(g)=\beta_g(f)$ for any
$g\in L_B$. Finally we define $\gamma^\prime: L_1(\mu_A)\times L_1(\mu_B)\to L_1(\mu_C)$  by 
$\gamma^\prime(f,g)={\gamma^f}^\prime(g)$ for any $f\in  L_1(\mu_A)$ and $g\in  L_1(\mu_B)$. It follows that 
$\gamma^\prime(f,g)=\widetilde{\gamma}(f,g)$ whenever  $f\in  L_A$ and $g\in  L_B$, so $\gamma^\prime$ is also bounded.

By \cite[253F, 253G]{Fre} there exists a unique bounded linear operator $\Omega: L_1(\lambda) \rightarrow L_1(\mu_C)$ such that $\Omega (f\ten g)= \gamma^\prime(f,g)$ for any $f\in  L_1(\mu_A)$ and $g\in  L_1(\mu_B)$.   The desired bounded linear  function is $\omega=\Omega\mid_T$.

To prove the uniqueness, let $\omega^\prime:T\to  L_1(\mu_C)$ be another bounded linear function that closes the above diagram.  Using \cite[Proposition 4.2]{FloLeu}  we extend it to the \textit{f}-algebra generated by $T$ in  $L_1(\lambda)$. Applying the B.L.T. theorem we get a bounded bilinear transformation $\Omega^\prime:  L_1(\lambda) \rightarrow L_1(\mu_C)$ such that 
$\Omega^{\prime}(f\ten g)= \gamma (f,g)$ for any $f\in  L_1(\mu_A)_u$ and $g\in  L_1(\mu_B)_u$. Following similar arguments as above one gets  $\Omega^{\prime}(f\ten g)= \gamma^\prime (f,g)$ for any $f\in  L_1(\mu_A)$ and $g\in  L_1(\mu_B)$, so $\Omega^\prime=\Omega$ and $\omega^\prime=\omega$. 

\end{proof}

\bfact

If $(A, s_A)$ and $(B, s_B)$  are two probability MV-algebras we defined a probability MV-algebra  $(T,s)$  and a bilinear  function $\beta: A\times B\rightarrow T$ such that  $(A, s_A)$ and $(B, s_B)$ are  $(T, s, \beta)$-independent.
Thorem \ref{teo:unpropInd} can be seen as a "universal property" of the product space, but it does not define $(T,s_T)$ up to isomorphism. If $(A,s_A)$ and $(B,s_B)$ are probability MV-algebras and both $(T,s_T)$ and $(V,s_V)$ satisfy the property from Theorem \ref{teo:unpropInd} then $T$ and $V$ have isomorphic group reducts. 
\efact

\bibliographystyle{elsarticle-num}

\begin{thebibliography}{24}

\bibitem{Birk}
Birkoff G., \textit{Lattice Theory}, AMS Coll. Publ. 25 3rd Ed. 1973.

\bibitem{BP}
Birkhoff G., Pierce R.S., \textit{Lattice-ordered Rings}, An. Acad. Brasil. Cienc. 28 (1956) 41-69.

\bibitem{Cha1}
Chang C.C., \textit{Algebraic analysis of many valued logics}, Trans. Amer. Math. Soc. 88 (1958) 467-490.

\bibitem{CDM}
Cignoli R., D'Ottaviano I.M.L., Mundici D., \textit{Algebraic foundation of many valued Reasoning}, Kluver Academc Publ (2000) Dordrecht.

\bibitem{Conrad}
Conrad P., \textit{The additive group of an f-ring}, Can. J. Math. 26(5) (1974) 1157-1168.

\bibitem{JVanR}
De Jonge E. Van Rooij, A.C.M., \textit{Introduction to Riesz Spaces}, Mathematical Centre Tracts 78, Amsterdam, 1977.

\bibitem{DiaLeu}
Diaconescu D., Leu\c stean I., \textit{The Riesz Hull of a Semisimple MV-algebra}, 

\bibitem{DiND}
Di Nola A., Dvurecenskij A., \textit{Product MV-algebras}, Multiple-Valued Logics 6 (2001) 193-215.


\bibitem{LDHandbook}
Di Nola A., Leu\c stean I., \textit{Handbook of Mathematical fuzzy logic - Volume 1}, volume 37 of Studies in Logic, chapter \L ukasiewicz logic and MV-algebras, College Publications, London 2011.

\bibitem{LeuRMV}
Di Nola A., Leustean I., \textit{\L ukasiewicz logic and Riesz Spaces}, Soft Comp. , accepted.


\bibitem{Fe}
Feller, W.,  \textit{An Introduction to Probability Theory and Its Applications}, Vol. II, Wiley, New York, 1971.



\bibitem{FloLeu}
Flondor P., Leu\c stean I., \textit{MV-algebras with operators (the commutative
and the non-commutative case)},  Discrete Mathematics 274 (2004) 41-76.

\bibitem{Fre}
Fremlin D.H., \textit{Measure Theory}, 1995. Avaliable at \textit{http://www.essex.ac.uk/maths/staff/fremlin/mt.htm}

\bibitem{Brunella}
Gerla B., \textit{Rational \L ukasiewicz logic and DMV-algebra}, Neural Networks World 11 (2001) 579-584.

\bibitem{Ha1}
Hausdorff  F., \textit{Summationmethoden und Momentfolgen I},  Math. Z. 9 (1921), 74 -109.

\bibitem{Ha2}
 Hausdorff F.,  \textit{Summationmethoden und Momentfolgen II},  Math. Z. 9 (1921), 280-299.

\bibitem{Kak}
Kakutani S., \textit{Concrete representation of abstract (L)-space and the mean ergodic theorem}, Annals of Mathematics 42 (1941) 523-537.

\bibitem{K}
Kroupa T.,  \textit{Every state on semisimple MV-algebra is integral}, Fuzzy Sets ad Systems,
157  (2006), 2771-2782.

\bibitem{KExt}
Kroupa T., \textit{Representation and extension of states on MV-algebras}, Arch. Math. Logic 45 (2006) 381-392.

\bibitem{LL1}
Lapenta S., Leu\c stean I., \textit{Towards Pierce-Birkhoff conjecture via MV-algebras}, submitted.

\bibitem{LeuStateC}
Leu\c stean I., \textit{Metric completion of MV-algebras: an approach to Stochastic Indepencence}, Journal of Logic and Computation 21(3) (2011) 493-508.

\bibitem{LeuSCRMV}
Leu\c stean I., \textit{State-complete Riesz MV-algebras and L-Measure Spaces}, Proceedings of IMPU, Springer (2012) 226-234.


\bibitem{Mir}
 Miranda E., de Cooman G., Quaeghebeur E., \textit{The Hausdorff moment problem under finite additivity}, Journal of Theoretical Probability 20(3) 2007 pp 663-693.

\bibitem{MonPMV}
Montagna F., \textit{An algebraic approach to Propositional Fuzzy Logic}, Journal of Logic, Language and Information 9 (2000) 91-124.

\bibitem{Mun1}
Mundici D., \textit{Interpretation of ACF*-algebras in \L ukasiewicz sentential calculus}, J. Funct. Anal. 65 (1986) 15-63.

\bibitem{Mu} 
Mundici D., \textit{Averaging the Truth-Value in \L ukasiewicz Logic}, Studia Logica 55 (1995), 113-127.

\bibitem{MunBook}
Mundici D., \textit{Advances in \L ukasiewicz calculus and MV-algebras}, Trends in Logic 35 Springer (2011).

\bibitem{P}
Panti G. \textit{ Invariant measures in free MV-algebras}, Communications in Algebra 36 (2008), 2849-2861.

\bibitem{BLT}
Reed, M.,   Simon B.,  \textit{Methods of Modern Mathematical Physics, Vol. 1: Functional Analysis},  Academic Press 1980.

\bibitem{MunRie}
Rie\v{c}an B., Mundici D., \textit{Probability on MV-algberas}, Chapter in the Handbook of Measure Theory, Elsevier Science 2002.



\end{thebibliography}

\end{document}